\documentclass[11pt,letterpaper,twoside,reqno]{amsart}

\usepackage{amsmath}
\usepackage{amssymb}
\usepackage{bm}
\usepackage{mathrsfs}
\usepackage[dvips]{graphicx}

\usepackage{color}

\usepackage{url}

\usepackage{supertabular}

\usepackage{alg}

\makeatletter
\def\algspacing{\alg@unmargin}
\makeatother

\newtheorem{thm}{Theorem}

\newtheorem{defn}[thm]{Definition}

\numberwithin{equation}{section}

\numberwithin{thm}{section}

\DeclareMathAlphabet{\mathsfsl}{OT1}{cmss}{m}{sl}

\newcommand{\term}{\emph}

\newcommand{\cnst}[1]{\mathrm{#1}}

	{\makebox{\phantom{}} \\ %
		\textsc{Input:}\begin{itemize}}%
	{\end{itemize}}

	{\textsc{Output:}\begin{itemize}}%
	{\end{itemize}}

	{\textsc{Procedure:}\begin{enumerate}}%
	{\end{enumerate}}

\renewcommand{\phi}{\varphi}

\newcommand{\eps}{\varepsilon}

\newcommand{\coll}[1]{\mathscr{#1}}

\newcommand{\abs}[1]{\left\vert {#1} \right\vert}

\newcommand{\argmin}{\operatorname*{arg\; min}}

\newcommand{\vct}[1]{\bm{#1}}
\newcommand{\mtx}[1]{\bm{#1}}

\newcommand{\adj}{*}
\newcommand{\psinv}{\dagger}

\newcommand{\supp}[1]{\operatorname{supp}(#1)}

\newcommand{\norm}[1]{\left\Vert {#1} \right\Vert}

\newcommand{\enorm}[1]{\norm{#1}_2}

\newcommand{\pnorm}[2]{\norm{#2}_{#1}}

\newcommand{\subjto}{\quad\text{subject to}\quad}

\newcommand{\bigO}{\mathrm{O}}

\newcommand{\Fee}{\mtx{\Phi}}

\newcommand{\restrict}[1]{\vert_{#1}}

\usepackage{subfigure}
\usepackage{wasysym}

\theoremstyle{remark}

\newcommand{\cosamp}{{\rm\sf CoSaMP}}

\begin{document}
\bibliographystyle{plain}
\title{Greedy Signal Recovery Review}

\author{Deanna Needell, 
Joel A. Tropp, 
Roman Vershynin
}
\thanks{DN is with the Mathematics Dept., Univ.~California at Davis, 1 Shields Ave., Davis, CA 95616.
E-mail: \url{dneedell@math.ucdavis.edu}.
RV is with the Mathematics Dept., Univ.~Michigan, Ann Arbor, MI 48109.
E-mail: \url{romanv@umich.edu}.
JAT is with Applied and Computational Mathematics, MC 217-50, California Inst.~Technology, Pasadena, CA 91125.
E-mail: \url{jtropp@acm.caltech.edu}.}

\date{July 2008}
\begin{abstract}
The two major approaches to sparse recovery are
$L_1$-minimization and greedy methods. Recently, Needell and Vershynin developed Regularized
Orthogonal Matching Pursuit (ROMP) that has bridged the gap between these two approaches. 
ROMP is the first stable greedy algorithm providing uniform guarantees. 

Even more recently, Needell and Tropp developed the stable greedy algorithm Compressive Sampling Matching
Pursuit (CoSaMP). CoSaMP provides uniform guarantees and improves upon
the stability bounds and RIC requirements of ROMP. CoSaMP offers rigorous bounds on computational cost
and storage. 
In many cases, the running time is just $O(N \log N)$,
where $N$ is the ambient dimension of the signal. This review summarizes these major advances.

\end{abstract}

\maketitle

\section{Introduction}
Sparse signals are those that contain much less information than their ambient dimension suggests.  The conventional signal compression scheme acquires the entire signal and then compresses it.  This methodology has been questioned for decades, and new approaches in compressed sensing have been developed to overcome this seemingly wasteful approach.  

Suppose $\vct{x}$ is a signal in $\mathbb{R}^N$, and define the $\ell_0$ quasi-norm\footnote{We consider real numbers here for simplicity, but similar results hold for the complex case.}
$$
\pnorm{0}{ \vct{x} } = \abs{ \supp{ \vct{x} }}= \abs{ \{ j : x_j \neq 0 \} }.
$$
When $\pnorm{0}{\vct{x}} \leq s$, we say that the signal $\vct{x}$ is \term{$s$-sparse}.  In practice, signals are not exactly sparse, but are rather close to sparse vectors. For example, compressible signals are those whose coefficients decay rapidly when sorted by magnitude.  We say that a signal $\vct{x}$ is \term{$p$-compressible} with magnitude $R$ if the sorted components of the signal decay at the rate
\begin{equation}\label{compressible}
\abs{x}_{(i)} \leq R \cdot i^{-1/p}
\qquad\text{for $i = 1, 2, 3, \dots$}.
\end{equation}

The sparse recovery problem is the reconstruction of such signals from a set of nonadaptive linear measurements.  The measurements are of the form $\Fee \vct{x}$ where $\Fee$ is some $m \times N$ measurement matrix.  Although in general this recovery is NP-Hard, work in compressed sensing has shown that for certain kinds of measurement matrices, recovery is possible when the number of measurements $m$ is nearly linear in the sparsity $s$, 
\begin{equation}   \label{N}                  
  m = s \log^{O(1)} (N).
\end{equation}
The survey~\cite{C06:CS} contains a discussion of these results.
 
The two major algorithmic approaches to sparse recovery are based on 
$L_1$-minimization and on greedy methods (Matching Pursuits). 
In this review we briefly describe these methods, as well as two new iterative methods that
provide the advantages of both approaches.  The first method, ROMP, is the first
stable algorithm to provide uniform guarantees. The second method, CoSaMP, improves
upon the results of ROMP, and provides rigorous computational bounds.

\section{Major approaches}
\subsection{$L_1$-minimization}

To recover the sparse signal $\vct{x}$ from its
measurements $\Fee \vct{x}$, one needs to find the the solution to the highly non-convex problem
\begin{equation}                    \tag{$L_0$}
  \min \|\vct{z}\|_0 \qquad \text{subject to} \qquad \Phi \vct{z} = \Phi \vct{x}.
\end{equation}

Donoho and his associates~\cite{DS89:Uncertainty-Principles} suggested
that for some measurement matrices $\Phi$, the generally NP-Hard problem $(L_0)$ should be equivalent 
to its convex relaxation: 
\begin{equation}                    \tag{$L_1$}
  \min \|\vct{z}\|_1 \qquad \text{subject to} \qquad \Phi \vct{z} = \Phi \vct{x},
\end{equation}
where $\|\vct{z}\|_1 = \sum_i |z_i|$ denotes the $\ell_1$-norm.  The convex problem $(L_1)$ can be solved using 
methods of linear programming. 

Clearly if the measurement matrix $\Fee$ is one-to-one on all $2s$-sparse vectors, then the $s$-sparse signal $\vct{x}$ will be recovered by solving $(L_0)$.  Cand\`es and Tao \cite{CT05:Decoding} proved that if $\Fee$ satisfies a stronger
condition then recovery is possible by solving the convex problem $(L_1)$.

\begin{defn}[Restricted Isometry Condition] 
  A measurement matrix $\Phi$ satisfies the 
  {\em Restricted Isometry Condition} (RIC)
  with parameters $(n, \delta)$ for $\delta \in (0,1)$
  if we have
  $$
  (1-\delta)\|\vct{v}\|_2 \leq \|\Phi \vct{v}\|_2 \leq (1+\delta)\|\vct{v}\|_2
  \qquad \text{for all $n$-sparse vectors}.
  $$
\end{defn}

When $\delta$ is small, the restricted isometry condition says that every set of $n$ columns of $\Fee$ is approximately an orthonormal system.  It has been shown (see \cite{RV08:sparse} and \cite{MPJ06:Uniform}) that random Gaussian, Bernoulli and partial Fourier matrices satisfy the restricted isometry condition with number of measurements as in \eqref{N}.  In the more practical case when $\vct{x}$ is not exactly sparse and corrupted with noise, we consider the mathematical program:
\begin{equation} \label{eqn:bp}
\min \pnorm{1}{ \vct{y} }
\quad\subjto\quad
\enorm{\Fee \vct{y} - \vct{u}}\leq \eps.
\end{equation}

This program approximately recovers the signal $\vct{x}$ even in the presence of noise.

\begin{thm}[Recovery under RIC \cite{CRT06:Stable}]  \label{stable RIC}
  Assume that the measurement matrix $\Phi$ satisfies the
  Restricted Isometry Condition with parameters $(3s, 0.2)$.
  Let $\Fee$ be a measurement matrix and let $\vct{u} = \Fee \vct{x} + \vct{e}$ be a noisy measurement vector where $\vct{x}\in\mathbb{R}^N$ is an arbitrary signal and $\enorm{\vct{e}}\leq \eps$. Then the program~\eqref{eqn:bp} produces an approximation $\vct{x}^\#$ that satisfies:
  $$
  \enorm{ \vct{x} - \vct{x}^\# } \leq \cnst{C} \left[ \frac{1}{\sqrt{s}} \pnorm{1}{\vct{x} - \vct{x}_s }
	+ \eps \right],
  $$
  where $\vct{x}_s$ denotes the $s$-sparse vector consisting of the $s$ largest coefficients in magnitude of $\vct{x}$. 
\end{thm}

In~\cite{Can08:Restricted-Isometry}, Cand\`es sharpened this theorem to work under the restricted isometry condition with parameters $(2s, \sqrt{2}-1)$.  This theorem also demonstrates that in the noiseless case, the $L_1$ approach provides exact reconstruction.  This was proved initially in~\cite{CT05:Decoding}.

\subsection{Orthogonal Matching Pursuit (OMP)}
An alternative approach to sparse recovery is via greedy algorithms. 
These methods find the support of the signal $\vct{x}$ iteratively, and reconstruct the signal using the pseudoinverse.

Orthogonal Matching Pursuit (OMP) is such an algorithm, analyzed by Gilbert and Tropp in~\cite{TG07:Signal-Recovery}.  Since we expect the columns of the measurement matrix $\Fee$ to be approximately orthonormal, $\Fee^*\Fee \vct{x}$ is locally a good approximation to $\vct{x}$.  OMP uses this idea to compute the support of a $s$-sparse signal $\vct{x}$.  First, the residual $\vct{r}$ is set to the measurement vector $\vct{u}$. At each iteration, the observation vector is set, $\vct{y} = \Phi^* \vct{r}$, and the coordinate of its largest coefficient in magnitude is added to the index set $I$. Then by solving a least squares problem, the residual is updated to remove this coordinate's contribution,
 $$
\vct{y} = \argmin_{\vct{z} \in \mathbb{R}^I} \|\vct{u} - \Phi \vct{z}\|_2; \qquad \vct{r} = \vct{u} - \Phi \vct{y}.
$$
Repeating this $s$ times yields an index set of $s$ coordinates corresponding to the support of the signal $\vct{x}$.
Tropp and Gilbert \cite{TG07:Signal-Recovery} showed that OMP recovers a sparse signal with high probability.

\begin{thm}[OMP Recovery~\cite{TG07:Signal-Recovery}] \label{OMP}
Let $\Fee$ be a $m \times N$ subgaussian matrix, and fix a $s$-sparse signal $\vct{x}\in\mathbb{R}^N$. Then OMP recovers (the support of) $\vct{x}$ from the measurements $\vct{u} = \Fee \vct{x}$ correctly with high probability, provided the number of measurements is $m \sim s \log N$.
\end{thm}

\subsection{Advantages and challenges of both approaches}
The $L_1$-mini\-mization method provides {\em uniform guarantees} for sparse recovery. 
Once the measurement matrix $\Phi$ satisfies the restricted isometry condition, 
this method works correctly for all sparse signals $\vct{x}$.  The method is also {\em stable}, so it works for non-sparse signals such as those which are compressible, as well as noisy signals.  However, the method is based on linear programming, and there is no strongly polynomial time algorithm in linear programming yet.

OMP on the other hand, is quite {\em fast} both provably and empirically.  The speed of OMP is a great advantage, but it lacks the strong guarantees that $L_1$ provides. Indeed, OMP works correctly for a {\em fixed} signal and measurement matrix with high probability, and so it must fail for some sparse signals and matrices~\cite{R08:Impossibility}.  It is also unknown whether OMP succeeds for compressible signals or on noisy measurements.

There has thus existed a gap between the approaches. The development of Regularized Orthogonal Matching Pursuit (ROMP) bridges this gap by providing a greedy algorithm with the same advantages as the $L_1$ method.  Compressive Sampling Matching Pursuit (CoSaMP) improves upon these results and provides rigorous runtime guarantees.  We now discuss these new algorithms.

\section{Regularized OMP}

Regularized Orthogonal Matching Pursuit (ROMP) is a new algorithm developed by Needell and Vershynin in~\cite{NV07:Uniform-Uncertainty,NV07:ROMP-Stable} for sparse recovery that performs correctly for all measurement matrices $\Phi$ satisfying the restricted isometry condition,
and for all sparse signals. Again, since the restricted 
isometry condition gaurantees every set of $s$ columns forms approximately an orthonormal system, every $s$ coordinates of the observation vector $\vct{y} = \Fee^* \vct{u}$ are in a loose sense good estimators of the corresponding $s$ coefficients of $\vct{x}$.  This notion suggests to use the $s$ largest coefficients of the observation vector $\vct{y}$ rather than only the largest, as in OMP.  We also include a regularization step to ensure that each coordinate carries close to an even share of the information.  The ROMP algorithm is described in Algorithm~\ref{alg:romp}.

\begin{algorithm}[thb]
\caption{ROMP Recovery Algorithm~\cite{NV07:Uniform-Uncertainty}}
	\label{alg:romp}
\centering \fbox{
\begin{minipage}{.9\textwidth} 
\vspace{4pt}
\algname{ROMP}{$\Fee$, $\vct{u}$, $s$}
\alginout{Measurement matrix $\Fee$, measurement vector $\vct{u}$, sparsity level $s$}
{Index set $I \subset \{1,\ldots,d\}$}
\vspace{8pt}\hrule\vspace{8pt}

\begin{algtab*}

\begin{description}
    \item[Initialize] Let the index set $I = \emptyset$ and the residual $r = u$.\\
      Repeat the following steps $s$ times or until $|I|\geq 2s$:
    \item[Identify] Choose a set $J$ of the $s$ biggest coordinates in magnitude 
      of the observation vector $\vct{y} = \Phi^*\vct{r}$, or all of its nonzero coordinates, 
      whichever set is smaller.
    \item[Regularize] Among all subsets $J_0 \subset J$ with comparable coordinates:
      $$
      |\vct{y}(i)| \leq 2|\vct{y}(j)| \quad \text{for all } i,j \in J_0,
      $$
      choose $J_0$ with the maximal energy $\|\vct{y}|_{J_0}\|_2$.
    \item[Update] Add the set $J_0$ to the index set: $I \leftarrow I \cup J_0$, 
      and update the residual:
      $$
      \vct{w} = \argmin_{\vct{z} \in \mathbb{R}^I} \|\vct{u} - \Phi \vct{z}\|_2; \qquad \vct{r} = \vct{u} - \Phi \vct{w}.
      $$
  \end{description}

\end{algtab*}
\vspace{-14pt}
\end{minipage}
}
\end{algorithm}

For measurement matrices that satisfy the RIC, ROMP stably recovers all $s$-sparse signals.  This is summarized in the following theorem from~\cite{NV07:ROMP-Stable}.

\begin{thm}[Recovery by ROMP~\cite{NV07:ROMP-Stable}]\label{T:stabsig}
  Assume a measurement matrix $\Phi$ satisfies the restricted isometry condition 
  with parameters $(8s, \epsilon)$ for $\epsilon = 0.01 / \sqrt{\log s}$. 
  Consider an arbitrary vector $\vct{x}$ in $\mathbb{R}^N$.
 Suppose that the measurement vector $\Phi \vct{x}$ becomes corrupted, 
  so we consider $\vct{u} = \Phi \vct{x} + \vct{e}$ where $\vct{e}$ is some error vector. 
  Then ROMP produces a good approximation $\hat{\vct{x}}$ to $\vct{x}$:
  \begin{equation}\label{boundsig}
    \|\hat{\vct{x}} - \vct{x}\|_2 
    \leq C \sqrt{\log s} \Big( \|\vct{e}\|_2 + \frac{\|\vct{x}-\vct{x}_s\|_1}{\sqrt{s}} \Big).
  \end{equation}  
\end{thm}

In the case where the signal $\vct{x}$ is exactly sparse without noise, this theorem guarantees exact reconstruction.  Note also that in the noisy case, ROMP needs no knowledge about the error vector $\vct{e}$ to approximate the signal.  In the special case where $\vct{x}$ is a compressible signal as in~\eqref{compressible}, the theorem provides the bound
\begin{equation}\label{boundcomp}
      \|\vct{x} - \hat{\vct{x}}\|_2 \le R' \frac{\sqrt{\log s}}{s^{p-1/2}} + C\sqrt{\log s}\|\vct{e}\|_2.
\end{equation}

ROMP thus provides the first greedy approach with these uniform and stable guarantees.  See~\cite{NV07:Uniform-Uncertainty,NV07:ROMP-Stable} for runtime analysis and empirical results.  The results are optimal up to the logarithmic factors, and the runtime is similar to that of OMP.  There is, however, room for improvement which leads us to the CoSaMP algorithm, developed by Needell and Tropp~\cite{NT08:Cosamp}.  CoSaMP provides optimal guarantees as well as an important implementation analysis.

\section{CoSaMP}

CoSaMP is an iterative recovery algorithm that provides the same guarantees as even the best optimization approaches.  As in the case of ROMP and the $L_1$ approach, CoSaMP recovers signals using measurement matrices that satisfy the RIC.  Thus as before, the observation vector $\vct{y} = \Fee^* \vct{u}$ serves as a good proxy for the signal $\vct{x}$.  Using the largest coordinates, an approximation to the signal is formed at each iteration.  After each new residual is formed, reflecting the missing portion of the signal, the measurements are updated.  This is repeated until all the recoverable portion of the signal is found.  (See~\cite{NT08:Cosamp} for halting criteria.)  The algorithm is described in Algorithm~\ref{alg:cosamp}.

\begin{algorithm}[thb]
\caption{\cosamp\ Recovery Algorithm~\cite{NT08:Cosamp}}
	\label{alg:cosamp}
\centering \fbox{
\begin{minipage}{.9\textwidth} 
\vspace{4pt}
\algname{\cosamp}{$\Fee$, $\vct{u}$, $s$}
\alginout{Sampling matrix $\Fee$, noisy sample vector $\vct{u}$, sparsity level $s$}
{An $s$-sparse approximation $\vct{a}$ of the target signal
}
\vspace{8pt}\hrule\vspace{8pt}

\begin{algtab*}
$\vct{a}^0 \leftarrow \vct{0}$, $\vct{v} \leftarrow \vct{u}$, $k \leftarrow 0$
	 	\hfill \{ Trivial initial approximation \} \\

\algrepeat
	$k \leftarrow k + 1$ \\
	$\vct{y} \leftarrow \Fee^\adj \vct{v}$
		\hfill \{ Form signal proxy \} \\
	$\Omega \leftarrow \supp{ \vct{y}_{2s} }$
		\hfill \{ Identify large components \} \\	
	$T \leftarrow \Omega \cup \supp{ \vct{a}^{k-1} }$
		\hfill \{ Merge supports \} \\
	
	$\vct{b}\restrict{T} \leftarrow \Fee_T^\psinv \vct{u}$
		\hfill \{ Signal estimation by least-squares \} \\		$\vct{b}\restrict{T^c} \leftarrow \vct{0}$ \\
	
	$\vct{a}^{k} \leftarrow \vct{b}_s$
		\hfill \{ Prune to obtain next approximation \} \\
	$\vct{v} \leftarrow \vct{u} - \Fee \vct{a}^{k}$
		\hfill \{ Update current samples \} \\
\alguntil{halting criterion {\it true}} \\
	

\end{algtab*}
\vspace{-14pt}
\end{minipage}}
\end{algorithm}

For a recovery algorithm to be used efficiently in practice, the least squares step must be analyzed carefully.  By the restricted isometry condition, the matrix $\Fee_T$ in the estimation step is very well-conditioned.  This suggests the use of an iterative method such as Richardson's iteration~\cite[Sec.~7.4]{Bjo96:Numerical-Methods} to apply the psuedoinverse $\Fee_T^\psinv = (\Fee_T^\adj \Fee_T)^{-1} \Fee_T^\adj$.  This method is analyzed in the context of CoSaMP and shown to provide an efficient means of acquiring the estimation~\cite{NT08:Cosamp}.  The following theorem from~\cite{NT08:Cosamp} summarizes the fundamentally optimal recovery guarantees and rigorous computational costs of CoSaMP.

\begin{thm}[Recovery by CoSaMP \cite{NT08:Cosamp}] \label{thm:Cosamp}
Suppose that $\Fee$ is an $m \times N$ measurement matrix satisfying the restricted isometry condition with parameters $(2s, \cnst{c})$.  Let $\vct{u} = \Fee \vct{x} + \vct{e}$ be a vector of samples of an arbitrary signal, contaminated with arbitrary noise.  For a given precision parameter $\eta$, the algorithm \cosamp\ produces an $s$-sparse approximation $\vct{a}$ that satisfies
$$
\enorm{ \vct{x} - \vct{a} } \leq
	\cnst{C} \cdot \max\left\{ \eta, 
	\frac{1}{\sqrt{s}} \pnorm{1}{\vct{x} - \vct{x}_{s/2}} + \enorm{ \vct{e} }.
	\right\}
$$
The running time is $\bigO( \coll{L} \cdot \log ( \enorm{\vct{x}} / \eta ) )$, where $\coll{L}$ bounds the cost of a matrix--vector multiply with $\Fee$ or $\Fee^\adj$.  Working storage is $\bigO(N)$.
\end{thm}

In~\cite{DT08:CoSaMP-TR} it is shown that only a fixed number of iterations is required to reduce the error to an optimal amount.  The report~\cite{DT08:CoSaMP-TR} also discusses variations on the CoSaMP algorithm that may improve results in practice.

\bibliography{cosamp}

\begin{thebibliography}{10}

\bibitem{Bjo96:Numerical-Methods}
{\AA}.~Bj{\"o}rck.
\newblock {\em Numerical Methods for Least Squares Problems}.
\newblock SIAM, Philadelphia, 1996.

\bibitem{C06:CS}
E.~Cand\`es.
\newblock Compressive sampling.
\newblock In {\em Proc. International Congress of Mathematics}, volume~3, pages
  1433--1452, Madrid, Spain, 2006.

\bibitem{CRT06:Stable}
E.~Cand\`es, J.~Romberg, and T.~Tao.
\newblock Stable signal recovery from incomplete and inaccurate measurements.
\newblock {\em Communications on Pure and Applied Mathematics},
  59(8):1207--1223, 2006.

\bibitem{CT05:Decoding}
E.~Cand\`es and T.~Tao.
\newblock Decoding by linear programming.
\newblock {\em IEEE Trans. Inform. Theory}, 51:4203--4215, 2005.

\bibitem{Can08:Restricted-Isometry}
E.~J. Cand\`es.
\newblock The restricted isometry property and its implications for compressed
  sensing.
\newblock {\em C. R. Math. Acad. Sci. Paris, Serie I}, 346:589--—592, 2008.

\bibitem{DS89:Uncertainty-Principles}
D.~L. Donoho and P.~B. Stark.
\newblock Uncertainty principles and signal recovery.
\newblock {\em SIAM J. Appl. Math.}, 49(3):906--931, June 1989.

\bibitem{MPJ06:Uniform}
S.~Mendelson, A.~Pajor, and N.~Tomczak-Jaegermann.
\newblock Uniform uncertainty principle for {B}ernoulli and subgaussian
  ensembles.
\newblock To appear, {\it Constr. Approx.}, 2009.

\bibitem{DT08:CoSaMP-TR}
D.~Needell and J.~A. Tropp.
\newblock {CoSaMP}: {I}terative signal recovery from incomplete and inaccurate
  samples.
\newblock {ACM} {T}echnical {R}eport 2008-01, California Institute of
  Technology, Pasadena, July 2008.

\bibitem{NT08:Cosamp}
D.~Needell and J.~A. Tropp.
\newblock {C}o{S}a{M}{P}: Iterative signal recovery from noisy samples.
\newblock {\em Appl. Comput. Harmon. Anal.}, 2008.
\newblock DOI: 10.1016/j.acha.2008.07.002.

\bibitem{NV07:ROMP-Stable}
D.~Needell and R.~Vershynin.
\newblock Signal recovery from incomplete and inaccurate measurements via
  regularized orthogonal matching pursuit.
\newblock Submitted for publication, October 2007.

\bibitem{NV07:Uniform-Uncertainty}
D.~Needell and R.~Vershynin.
\newblock Uniform uncertainty principle and signal recovery via regularized
  orthogonal matching pursuit.
\newblock DOI: 10.1007/s10208-008-9031-3, 2007.

\bibitem{R08:Impossibility}
H.~Rauhut.
\newblock On the impossibility of uniform sparse reconstruction using greedy
  methods.
\newblock To appear, {\it Sampl. Theory Signal Image Process.}, 2008.

\bibitem{RV08:sparse}
M.~Rudelson and R.~Vershynin.
\newblock On sparse reconstruction from {F}ourier and {G}aussian measurements.
\newblock {\em Comm. Pure Appl. Math.}, 61:1025--1045, 2008.

\bibitem{TG07:Signal-Recovery}
J.~A. Tropp and A.~C. Gilbert.
\newblock Signal recovery from random measurements via orthogonal matching
  pursuit.
\newblock {\em IEEE Trans. Info. Theory}, 53(12):4655--4666, 2007.

\end{thebibliography}

\end{document}